\documentclass{article}

\usepackage{amsmath}
\usepackage{amssymb}

\newtheorem{theorem}{Theorem}
\newtheorem{lemma}{Lemma}
\newtheorem{corollary}{Corollary}
\newtheorem{proposition}{Proposition}
\newtheorem{definition}{Definition}

\newtheorem{remark}{Remark}

\begin{document}

\title{Factorization of quasi-variational relations systems}
\author{Daniela Inoan}
\date{}

\maketitle

\begin{abstract}
Variational relation problems allow a general approach for variational inequalities, equilibrium problems, optimization problems, variational inclusions.
 In this paper we consider a system of quasi-variational relations and determine some conditions in which the solvability of the independent problems imply the existence of a solution for the system. We particularize then the result for a system of variational inequalities and for a constrained Nash equilibrium problem.
\end{abstract}

\emph{Accepted for publication in Acta Mathematica Vietnamica}

\section{Introduction}

In \cite{Lu2008}, D.T. Luc introduced a general model, called a variational relation problem, showing that it is a unifying approach for several equilibrium problems, optimization problems, variational or differential inclusions problems. The study of variational relations was continued in different papers, for instance  M. Balaj and L.J. Lin \cite{Ba-Li2011}, M. Balaj and D.T. Luc \cite{Ba-Lu2010}, P.Q. Khanh and D.T. Luc \cite{Kh-Lu2008}, L.J. Lin and Q.H. Ansari \cite{Li-An2010}. This general model can be considered, in nature, equivalent to a variational inclusion problem, see for instance N.X. Hai and P.Q. Khanh \cite{Ha-Kh2007new}.

Various kinds of variational inequality systems, systems of quasi-variational inclusions or systems of quasi-equilibrium problems were studied in numerous papers, see for instance N.X. Hai and P.Q. Khanh \cite{Ha-Kh2007},  D. Inoan \cite{In2009}, G. Kassay, J. Kolumb\' an and Z. P\'ales \cite{Ka-Ko-Pa2002}, J.L. Lin and C.I. Tu \cite{Li-Tu2008}, S. Plubtieng and K. Sombut \cite{Pl-So2010}  and the references therein. A natural unifying model for several problems of this type is a system of quasi-variational relations.

Let $I=\{1, \dots,n\}$ be an index set. For each $i \in I$, let $X_i$ be a nonempty  subset of a  real Hausdorff locally convex  space and $X=\prod_{i\in I} X_i$.
Let $S_i, Q_i : X \to 2^{X_i}$ be set-valued maps with nonempty values.
Consider also $R_i(x,y_i) $ to be a relation linking $x \in X$ and $y_i \in X_i$.

The following system of quasi-variational relations was introduced in \cite{Li-An2010} by L.J. Lin and Q.H. Ansari:
\begin{eqnarray*}
(SQVR) & \mbox{Find}\ \bar x=(\bar x_1, \dots, \bar x_n) \in X \ \mbox{such that for each}\ i\in I,\\ & \bar x_i \in S_i(\bar x)\ \mbox{and}\
 R_i(\bar x, y_i) \ \mbox{holds for all} \ y_i \in Q_i(\bar x).
\end{eqnarray*}
In \cite{Li-An2010} some existence results for a solution of $(SQVR)$ are established, using a maximal element theorem for a family of set-valued maps.

In this paper we want to investigate in what conditions the solvability of some independent variational relation problems implies the existence of a solution for the system. This  question was addressed in the case of variational inequality systems in \cite{Ka-Ko-Pa2002}. Following the ideas from \cite{Ka-Ko-Pa2002}, let $i \in I$ be a fixed index and let $x_j \in X_j$ be fixed, for $j\neq i$. We formulate the corresponding $i$th subproblem:
\begin{eqnarray*}
(QVR_i) & \mbox{Find}\ \bar x_i \in X_i \ \mbox{such that}\   \bar x_i \in S_i(x_1, \dots,\bar x_i, \dots, x_n)\ \\ & \mbox{and}\
 R_i(x_1, \dots,\bar x_i, \dots, x_n, y_i) \ \mbox{holds for all} \ y_i \in Q_i(x_1, \dots,\bar x_i, \dots, x_n).
\end{eqnarray*}

Using these independent problems we will give in Theorem \ref{teo:main} of the next Section a factorization result for the system (SQVR). Then we will particularize it for a system of Minty variational inequalities and for a constrained Nash
equilibrium problem.

We recall that (see for instance \cite{Hu-Pa1997}) if $U$ and $V$ are topological spaces, a set-valued mapping with nonempty values $T:U \to 2^V$ is said to be \emph{upper semicontinuous} (u.s.c.) if for every closed set $C \subset V$, the set $\{ u \in U \ |\ T(u) \cap C \neq \emptyset \}$ is closed. $T$ is said to be \emph{lower semicontinuous} (l.s.c.) if for every closed set $C \subset V$, the set $\{ u \in U \ |\ T(u) \subset C  \}$ is closed.

\begin{lemma}
Let $U$ and $V$ be Hausdorff  topological spaces. The mapping $T:U \to 2^V$ is l.s.c. on U if and only if, for any $u \in U$, any net $(u_\alpha)_{\alpha \in A }$ converging to $u$, for every $v \in T(u)$ there exists $v_\alpha \in T(u_\alpha)$ such that the net $(v_\alpha)_{\alpha \in A }$ converges to $v$.
\end{lemma}

\section{A factorization result}

To prove the main result of the paper, we will use the  Kakutani-Fan-Glicksberg fixed point theorem:
\begin{theorem}\label{teo:KFG}\cite{Ze1986}
Let $K$ be a nonempty compact convex subset of a locally
convex  space. If $T:K \to 2^K$ is upper semicontinuous and $T(x)$ is nonempty convex and closed for any
$x \in K$, then there exists a fixed point $\bar x$ of $T$ in $K$, $\bar x \in T(\bar x)$.
\end{theorem}
In what follows, we will denote the convex hull of a set $A$ by $\mbox{co}A$. For $i\in I$ and $x_j \in X_j$ fixed, $i\neq j$, denote the set of fixed points for the partial function $S_i(x_1, \dots, \cdot, \dots, x_n)$ by
\[
\textrm{Fix}\, S_i(x_1, \dots, \cdot, \dots, x_n)=\{ z_i \in X_i \ |\ z_i \in S_i(x_1, \dots,x_{i-1},z_i, x_{i+1}, \dots, x_n) \}.
\]
\begin{theorem}\label{teo:main}
Suppose that, for each $i \in I$,  $X_i$ is nonempty compact convex and that for any fixed $x_j \in X_j$, $j \neq i$ we have:

(H1) the problem $(QVR_i)$ admits a solution;

(H2) the set $\textrm{Fix}\,S_i(x_1, \dots,\cdot,  \dots, x_n)$ is convex;

(H3) if the relation $R_i(x_1, \dots,z_{i}^k, \dots, x_n, t_{i}^k)$ holds for  some $z_{i}^k \in X_i$,  for every $t_{i}^k \in Q_i(x_1, \dots,z_{i}^k, \dots, x_n)$, $k \in \{1,\dots, l\}$, then $R_i(x_1, \dots,z_{i}, \dots, x_n, \theta_{i})$ holds for every $z_i \in \mbox{co} \{ z_{i}^1, \dots,z_{i}^l\}$ and $\theta_{i} \in Q_i(x_1, \dots,z_{i}, \dots, x_n)$;

(H4) for every net $(x^\alpha)_{\alpha \in A}$  with $x^\alpha=(x^\alpha_1, \dots, x^\alpha_n) \in X$, such that $x^\alpha \to x^0=(x^0_1, \dots, x^0_n) \in X$, if $x^\alpha_i \in S_i (x^\alpha)$ for every $\alpha \in A$, then
$x^0_i \in S_i (x^0)$;

(H5) for every net $(x^\alpha)_{\alpha \in A} \subset X$  such that $x^\alpha \to x^0=(x^0_1, \dots, x^0_n)\in X$, if $R_i(x^\alpha,t_i)$ holds for every $t_i \in Q_i(x^\alpha)$, for every $\alpha \in A$, then $R_i(x^0,\theta_i)$ holds for every $\theta_i \in Q_i(x^0)$.

Then the problem $(SQVR)$ admits a solution.
\end{theorem}
\textbf{Proof:}
For a fixed index $i \in I$, let the set valued mapping $T_i : X \to 2^{X_i} $ be defined by:
\begin{eqnarray*}
& T_i(x_1, \dots, x_n)=  \mbox{Fix}\,S_i(x_1, \dots,\cdot,  \dots, x_n)\cap \\
&  \cap \{ z \in X_i \ |\ R_i(x_1, \dots,z, \dots, x_n, t_i) \
\mbox{holds for every } t_i \in Q_i(x_1, \dots,z, \dots, x_n)\}.
\end{eqnarray*}
Define $T : X \to 2^{X} $ by
\[
T(x_1, \dots, x_n)=T_1(x_1, \dots, x_n) \times \dots \times T_n(x_1, \dots, x_n)
\]
It is clear that any fixed point of this set-valued mapping is a solution of the problem (SQVR).

From $(H1)$ it follows that, for any $i \in I$ and $x \in X$, $T_i(x) \neq \emptyset$.

Also, from $(H2)$ and $(H3)$ we can prove that $T_i(x)$ is convex.

The graph of the mapping $T_i$ is closed.  Indeed, let $(x^\alpha_1, \dots, x^\alpha_n) \in X$ and $z_i^\alpha \in T_i(x^\alpha)$ be such that $(x^\alpha_1, \dots, x^\alpha_n) \to (x_1^0, \dots,x_n^0)$ and $z^\alpha_i \to z^0_i$. This means that $z^\alpha_i \in S_i (x^\alpha_1, \dots, z^\alpha_i, \dots, x^\alpha_n)$ and $R_i(x^\alpha_1, \dots, z^\alpha_i, \dots, x^\alpha_n,t)$ holds for all $t\in Q_i (x^\alpha_1, \dots, z^\alpha_i, \dots, x^\alpha_n)$. Using  the hypotheses (H4) and (H5) for the net $(x_1^\alpha, \dots, z_i^\alpha, \dots, x_n^\alpha) \to (x_1^0, \dots, z_i^0, \dots, x_n ^0)$ we get that $z^0_i \in S_i (x^0_1, \dots, z^0_i, \dots, x^0_n)$ and $R_i(x^0_1, \dots, z^0_i, \dots, x^0_n,\theta)$ holds for every $\theta \in Q_i (x^0_1, \dots, z^0_i, \dots, x^0_n)$, that is $(x^0_1, \dots,x_n^0,z^0_i) \in \textrm{Graph}\,T_i$.

We also have that  $T_i(x)$ is closed, for every $x \in X$.

Since the values of $T_i$ are closed, contained in a compact set and the graph is closed, then  $T_i$ is upper semicontinuous (see for instance \cite{Hu-Pa1997}).  Then also $T$ is upper semicontinuous and has nonempty convex closed values. Applying now the Kakutani-Fan-Glicksberg fixed point theorem we get the existence of a fixed point for $T$, which is also a solution of (SQVR).

The existence result obtained in \cite{Li-An2010} for the system $(SQVR)$ is proved in different conditions. These are imposed mostly on the relation $R_i$ and mappings $S_i$, $Q_i$ with $x$ as a global variable. The partial problems $(QVR_i)$ may have  solutions even in the absence of such conditions.

\begin{remark}
Hypothesis (ii) of Theorem 1.3 in \cite{Li-An2010} is a classic property of  closedness of relations (see also \cite{Lu2008}): For a fixed $y_i \in X_i$, $R_i(\cdot,y_i)$ is said to be closed in the first variable if for every net $(x^\alpha)_{\alpha \in A} \subset X$ that converges to $x^0\in X$, whenever $R_i(x^\alpha,y_i)$ holds for all $\alpha\in A$, then $R_i(x^0,y_i)$ holds too.

Our hypothesis $(H5)$ is distinct from this property, as can be seen from the following example. Let $n=2$, $X_1=X_2=\mathbf{R}$ and define $\varphi :\mathbf{R} \times \mathbf{R} \to \mathbf{R} $ by

\[
\varphi(x)= \left\{ \begin{array}{ll} x_1+x_1x_2+x_2, & \textrm{if}\ x_2 \neq 0, \\ -1, & \textrm{if}\ x_1 \neq 0, x_2=0 \\
  0, & \textrm{if}\ x_1=x_2=0, \end{array}\right.
\]
for any $x=(x_1,x_2)$.

$R_1(x,y_1)$ holds if and only if $\varphi(x_1,x_2) \leq \varphi(y_1, x_2)$.

Define also $Q_1 :\mathbf{R} \times \mathbf{R} \to 2^ \mathbf{R} $ by
\[
Q_1(x)= \left\{ \begin{array}{ll} D, & \textrm{if}\ x_2 \neq 0, \\ (0, +\infty), & \textrm{if}\ x_1 \neq 0, x_2=0 \\
  \{0\}, & \textrm{if}\ x_1=x_2=0, \end{array} \right.
\]
with $D$ an arbitrary subset of $\mathbf{R}$.
The relation defined above is not closed in the first variable for every $y_1$ (for instance fix $y_1=1$ and take the sequence $(\frac{1}{n}, \frac{1}{n}) \to (0,0)$) but hypothesis $(H5)$ is satisfied.

Convexity of relations can also be defined (see for instance \cite{Ba-Lu2010}); $R_i$ is said to be convex if whenever $R_i(x^k, y_i^k)$ holds for $x^k \in X$, $y_i^k \in X_i$, $k=1,2$, the relation $R_i(\lambda x^1+(1-\lambda)x^2, \lambda y_i^1+(1-\lambda)y_i^2)$ holds for any $\lambda \in [0,1]$. It can be proved that in the preceding example, relation $R_1$ is not convex, but (together with $Q_1$) it satisfies hypothesis $(H3)$.
\end{remark}

Knaster-Kuratowski-Mazurkiewicz mappings  play an important role in nonlinear analysis.
In a natural way, one can define also  KKM relations  (see  \cite{Lu2008}).
\begin{definition}
Let $i \in I$ be a fixed index and $x_j \in X_j$ fixed, for $j\neq i$. We say that $R_i(x_1, \dots, \cdot, \dots, x_n, \cdot)$ is KKM if
 for any $x_i^k \in X_i$, $k \in \{1, \dots, l\}$ and $\tilde x_i \in \mbox{co} \{ x_i^1, \dots, x_i^l \}$ there exists an index $j$ such that $R_i(x_1, \dots, \tilde x_i, \dots, x_n, x_i^j)$ holds.
\end{definition}
\begin{definition}
Let $i \in I$ be a fixed index, $x_j \in X_j$ fixed, for $j\neq i$ and $t \in X_i$ fixed. We say that the partial relation $R_i(x_1, \dots, \cdot, \dots, x_n, t)$ is closed if
for every net $(x_i^\alpha)_{\alpha \in A} \subset X_i$ with $x_i^\alpha \to x_i^0$, if $R_i(x_1, \dots, x_i^\alpha, \dots, x_n, t)$ holds for every $\alpha\in A$, then  $R_i(x_1, \dots, x_i^0, \dots, x_n, t)$ also holds.
\end{definition}
The next proposition gives an existence result for the partial problem $(QVR_i)$.

\begin{proposition}
Let $i \in I$ be a fixed index and $x_j \in X_j$ fixed, for $j\neq i$. Suppose that:

(H6) $X_i$ is nonempty convex compact;

(H7) for every $z_i \in X_i$, $\mbox{co} Q_i(x_1, \dots,z_i,\dots,x_n) \subset S_i(x_1, \dots,z_i,\dots,x_n)$;

(H8) the set $ \{ z \in X_i \ |\ y_i \in Q_i(x_1, \dots, z, \dots, x_n) \}$ is open in $X_i$, for every $y_i \in X_i$;

(H9) the set $\mbox{Fix}\,S_i(x_1, \dots,\cdot,  \dots, x_n)$ is closed;

(H10) for any $t \in X_i$ fixed, the relation $R_i(x_1, \dots, \cdot, \dots, x_n, t)$ is closed;

(H11) the relation $R_i(x_1, \dots, \cdot, \dots, x_n, \cdot)$ is KKM.

Then the partial problem $(QVR_i)$ admits a solution.
\end{proposition}
\textbf{Proof:}
It follows as  a direct consequence of Corollary 3.1 from \cite{Lu2008}, applied to the partial functions $S_i(x_1, \dots, \cdot, \dots, x_n)$,
$Q_i(x_1, \dots, \cdot, \dots, x_n)$ and to the partial relation $R_i(x_1, \dots, \cdot, \dots, x_n, \cdot)$.

\begin{remark}
a) Obviously, $(H4)$ implies $(H9)$;

b) As it is mentioned in \cite{Lu2008}, the condition of $X_i$ being compact can be replaced by a suitable coercivity condition.
\end{remark}

\section{Applications}
I. Let $X_i$ be nonempty linear subspaces of some real Hausdorff locally convex   spaces. Let $Z_i$ be real Hausdorff topological vector spaces, $\langle \cdot,\cdot \rangle _i$ continuous bilinear functions defined on $Z_i \times X_i$ and $F_i : X \to 2^{Z_i}$  set-valued mappings.

The \emph{Minty variational inequality problem} concerning the system of set-valued functions $F_i$ is the following (see \cite{Ka-Ko-Pa2002}):
\begin{eqnarray*}
(MVIP) & \mbox{Find}\ \bar x=(\bar x_1, \dots, \bar x_n) \in X \ \mbox{such that for each}\ i\in I,\\ &
 \mbox{for all} \ y_i \in X_i, \ \mbox{for all}\ f \in F_i(\bar x_1, \dots, y_i, \dots,\bar x_n)\ \langle f, y_i- \bar x_i \rangle_i \geq 0.
\end{eqnarray*}
To convert $(MVIP)$ to a special case of $(SQVR)$, we take the set-valued maps $S_i(x)=Q_i(x)=X_i$, for every $x=(x_1, \dots, x_n)\in X$. The relation $R_i(x, y_i)$ holds if and only if
$
\langle f, y_i-x_i \rangle_i \geq 0, \ \mbox{for every} \ f \in F_i(x_1, \dots,y_i,\dots, x_n)
$.

The partial problem, for $i \in I$ and $x_j \in X_j$ fixed,  $j\neq i$, is in this case
\begin{eqnarray*}
 (MVIP)_i & \mbox{Find}\ \bar x_i \in X_i \ \mbox{such that}\\  &   \mbox{ for all} \ y_i \in X_i, \mbox{for all}\ f \in F_i( x_1, \dots, y_i, \dots, x_n)\ \langle f, y_i-\bar x_i \rangle_i \geq 0.
\end{eqnarray*}
As a direct consequence of Theorem \ref{teo:main} we can obtain Theorem 1 from \cite{Ka-Ko-Pa2002}:
\begin{corollary}
Suppose that, for each $i \in I$,  $X_i$ is compact and that:

(a) for any fixed $x_j \in X_j$, $j \neq i$, the problem $(MVIP)_i$ admits a solution;

(b) for any fixed $y_i \in X_i$ the set-valued function
\[
(x_1,\dots,x_{i-1},x_{i+1}, \dots, x_n) \mapsto F_i(x_1,\dots,x_{i-1},y_i,x_{i+1}, \dots, x_n)
 \]
 is lower semicontinuous on $X_1 \times \dots \times X_{i-1}\times X_{i+1} \times \dots \times X_n$.

Then the problem $(MVIP)$ admits a solution.
\end{corollary}
\textbf{Proof:}
Hypotheses (H1), (H2) and (H4) are obviously satisfied.

To check (H3), let $z_{i}^k \in X_i$, $k \in \{1,\dots, l\}$ such that $R_i(x_1, \dots,z_{i}^k, \dots, x_n,t_{i}^k)$ holds for every $t_{i}^k \in Q_i(x_1, \dots,z_{i}^k, \dots,x_n)=X_i$, that is $\langle f, t_{i}^k-z_{i}^k \rangle_i \geq 0$ for every $f\in F_i(x_1, \dots,t_{i}^k, \dots,x_n)$.
Let $z_i=\sum_{k=1}^l \lambda^k z_{i}^k $ belong to $\mbox{co} \{ z_{i}^1, \dots, z_{i}^l \}$ and consider $\theta_i \in Q_i(x_1, \dots,z_{i}, \dots,x_n)=X_i$. For every $f\in F_i(x_1, \dots,\theta_{i}, \dots,x_n)$ we have $\langle f, \theta_{i}-z_{i}^k \rangle_i \geq 0$. Multiplying each inequality by $\lambda^k$ and summing by $k \in \{1, \dots, l\}$ it follows that $\langle f, \theta_{i}-z_{i} \rangle_i \geq 0$. So $R_i(x_1,\dots,z_i,\dots,x_n, \theta_i)$ holds.

Finally, (H5) can be obtained from condition (b). Indeed, consider a net $x^\alpha =(x^\alpha_1, \dots, x^\alpha_n)$ in $X$, $\alpha \in A$, converging to $x^0=(x^0_1, \dots,x^0_n)$ and such that for every $t_i \in Q_i(x^\alpha)=X_i$, for every $f \in F_i(x^\alpha_1,\dots,t_i, \dots, x^\alpha_n)$ we have $\langle f,t_i-x^\alpha_i \rangle_i \geq 0$. Let $\theta_i \in Q_i(x^0)=X_i$ and $g \in F_i(x^0_1,\dots,\theta_i, \dots, x^0_n)$. Due to the lower semicontinuity imposed in (b), there exist $g^\alpha \in F_i(x^\alpha_1,\dots,\theta_i, \dots, x^\alpha_n)$ such that $g^\alpha \to g$. Also we have from above that $\langle g^\alpha,\theta_i-x^\alpha_i \rangle_i \geq 0$ for all $\alpha \in A$. So using the continuity of $\langle \cdot, \cdot \rangle_i$ yields $\langle g,\theta_i-x^0_i \rangle_i \geq 0$, that is $R_i(x^0,\theta_i)$ holds.

II. Another particular case of the system (SQVR) is the \emph{constrained Nash equilibrium problem} (see for instance \cite{Li-An2010}).
For each $i \in I$, let $\varphi_i: X \to \mathbf{R}$ and consider the problem
\begin{eqnarray*}
(CNEP) & \mbox{Find}\ \bar x=(\bar x_1, \dots, \bar x_n) \in X \ \mbox{such that for each}\ i\in I,\\ & \bar x_i \in S_i(\bar x)\ \mbox{and}\
\varphi_i(\bar x_1, \dots,\bar x_n) \leq \varphi_i(\bar x_1, \dots, \bar x_{i-1}, y_i, \bar x_{i+1}, \dots, \bar x_n) \\
 & \mbox{for all} \ y_i \in Q_i(\bar x).
\end{eqnarray*}
In this case,  the relation $R_i(x,y_i)$ holds if and only if
\[
\varphi_i(x_1, \dots,x_n) \leq \varphi_i(x_1, \dots, x_{i-1}, y_i, x_{i+1}, \dots, x_n).
\]
The partial problem, for $i \in I$ and $x_j \in X_j$ fixed,  $j\neq i$ is then
\begin{eqnarray*}
 (CNEP)_i & \mbox{Find}\ \bar x_i \in X_i \ \mbox{such that}\   \bar x_i \in S_i(x_1, \dots,\bar x_i, \dots, x_n)\ \mbox{and}\\ & 
\varphi_i(x_1, \dots,  x_{i-1}, \bar x_i, x_{i+1}, \dots,  x_n) \leq \varphi_i(x_1, \dots,  x_{i-1}, y_i, x_{i+1}, \dots,  x_n) \\ &  \mbox{ for all} \ y_i \in Q_i(x_1, \dots,\bar x_i, \dots, x_n).
\end{eqnarray*}
From Theorem \ref{teo:main} we have:
\begin{corollary}
Suppose that, for each $i \in I$,  $X_i$ is nonempty compact convex, hypothesis (H2), (H4) hold and:

(a) for any fixed $x_j \in X_j$, $j \neq i$, the problem $(CNEP)_i$ admits a solution;

(b) the application $\varphi_i(x_1, \dots, x_{i-1}, \cdot, x_{i+1}, \dots, x_n)$ is quasiconvex and $Q_i$ does not depend of the i-th variable;

(c) $\varphi_i$ is continuous and $Q_i$ is lower semi-continuous.

Then the problem $(CNEP)$ admits a solution.
\end{corollary}
\textbf{Proof:}
Let $z_{i}^k \in X_i$, $k \in \{1,\dots, l\}$ be such that for all $t_{i}^k \in Q_i(x_1, \dots,z_{i}^k,\dots,x_n)$ the inequality $\varphi_i(x_1,\dots,z_{i}^k, \dots,x_n) \leq \varphi_i(x_1,\dots,t_{i}^k, \dots,x_n)$ holds. Let $z_i=\sum_{k=1}^l \lambda^k z_{i}^k$, with $\lambda^k\geq 0$, $\sum_{k=1}^l \lambda^k=1$ and $\theta_i \in Q_i(x_1, \dots,z_{i}, \dots,x_n)$. Since we have $Q_i(x_1, \dots,z_{i},\dots,x_n)=Q_i(x_1, \dots,z_{i}^k,\dots,x_n)$ for each $k$, it follows that $\varphi_i(x_1,\dots,z_{i}^k, \dots,x_n) \leq \varphi_i(x_1,\dots,\theta_{i}, \dots,x_n)$ and from the quasiconvexity of the partial function it follows that $\varphi_i(x_1,\dots,z_{i}, \dots,x_n) \leq \varphi_i(x_1,\dots,\theta_{i}, \dots,x_n)$. This insures (H3).

To check (H5), let $x^\alpha =(x^\alpha_1, \dots, x^\alpha_n)$ in $X$ be a net converging to $x^0=(x^0_1, \dots,x^0_n)$ and such that $R_i(x^\alpha,t_i)$ holds for every $t_i \in Q_i(x^\alpha)$, for every $\alpha \in A$. Let $\theta_i \in Q_i (x^0)$. From the lower semicontinuity of $Q_i$ it follows that there exists $t^\alpha_i \in Q_i(x^\alpha)$ such that $t^\alpha_i \to \theta_i$. Since
$t^\alpha_i \in Q_i(x^\alpha)$ we have $\varphi_i(x_1^\alpha,\dots,x_n^\alpha) \leq \varphi_i(x_1^\alpha,\dots,t^\alpha_i,\dots, x_n^\alpha)$.

But $\varphi_i$ is continuous, so also $\varphi_i(x_1^0,\dots,x_n^0) \leq \varphi_i(x_1^0,\dots,\theta_i ,\dots, x_n^0)$.

\begin{remark}
The classic case when the set-valued mapping $Q_i$ is given as
\[
Q_i(x_1, \dots, x_n)=\{ t \in X_i \ |\ g_i(x_1, \dots, t, \dots, x_n) \leq 0 \},
\]
 where $g_i:X \to \mathbf{R}$ is a given function, is an example for $Q_i$ independent with respect to the i-th variable.

For $X_i=\mathbf{R}$ ($i \in I$) and $g_i(x_1, \dots,x_n)= x_i- f_i(x_1, \dots, x_{i-1}, x_{i+1}, \dots, x_n)$, with $f: \mathbf{R}^{n-1} \to \mathbf{R}$ a continuous function, the set-valued mapping $Q_i$ becomes $Q_i(x_1, \dots, x_n)=(-\infty, f_i(x_1, \dots, x_{i-1}, x_{i+1}, \dots, x_n) ]$, and is also lower semicontinuous.
\end{remark}

Recently, new theorems on maximal elements, coincidence or fixed points, nonempty intersections of multivalued mappings were obtained in more general settings - for instance in generalized finitely continuous (GFC) topological spaces, without a linear structure (see N.X. Hai, P.Q. Khanh, N.H. Quan \cite{Ha-Kh-Qu2009}, P.Q. Khanh, N.H. Quan \cite{Kh-Qu2010} and the references therein). These results give the possibility of existence investigation with more relaxed assumptions, especially those related to convexity. The study of quasi-variational relations systems could be further extended in such general spaces, using these new tools to establish existence results.

\textbf{Acknowledgements} The author
is grateful to the referee for the valuable remarks which helped to improve the paper.

\begin{flushright}
Daniela Ioana Inoan \\
Technical University of Cluj-Napoca\\
 Department of Mathematics,
 Cluj-Napoca, Romania \\
Daniela.Inoan@math.utcluj.ro
\end{flushright}

\end{document}